\documentclass[12pt]{article}
\usepackage[UKenglish]{babel}
\usepackage[T1]{fontenc}
\usepackage{lmodern,amsmath,amsthm,amsfonts,amssymb,graphicx,float,microtype,thmtools,underscore,mathtools,thm-restate}
\usepackage[shortlabels]{enumitem}
\setlist[itemize]{topsep=0ex,itemsep=0ex,parsep=0ex}
\setlist[enumerate]{topsep=0ex,itemsep=0ex,parsep=0ex}
\usepackage[dvipsnames,svgnames,table]{xcolor}
\usepackage{todonotes}
\usepackage{xstring}
\usepackage{stackengine}
\stackMath
\usepackage{mathtools}
\usepackage{tikz}
\usetikzlibrary{arrows.meta}
\usetikzlibrary{calc}
\usepackage{tikzmacros}
\usepackage[unicode=true]{hyperref}
\hypersetup{ 
colorlinks,
linkcolor={blue!60!black},
citecolor={black},
urlcolor={blue!60!black},
pdftitle={Quasi-Isometry: Separating Multiplicative Distortion from Additive Distortion}}
\usepackage[capitalise, compress, nameinlink, noabbrev]{cleveref}
\crefname{lem}{Lemma}{Lemmas}
\crefname{thm}{Theorem}{Theorems}
\crefname{ques}{Question}{Theorems}
\crefname{cor}{Corollary}{Corollaries}
\crefname{enumi}{Item}{Items}
\crefformat{equation}{(#2#1#3)}
\Crefformat{equation}{Equation #2(#1)#3}
\crefformat{enumi}{#2#1#3}
\Crefformat{enumi}{Item (#2#1#3)}
\newcommand{\defn}[1]{\textcolor{Maroon}{\emph{#1}}}
\usepackage[longnamesfirst,numbers,sort&compress]{natbib}
\makeatletter
\def\NAT@spacechar{~}
\makeatother
\usepackage[margin=28mm]{geometry}
\renewcommand{\baselinestretch}{1.1}
\allowdisplaybreaks
\renewcommand{\epsilon}{\varepsilon}

\renewcommand{\geq}{\geqslant}
\renewcommand{\leq}{\leqslant}
\DeclareMathOperator{\dist}{dist}

\newlength\mystringlen
\newcommand{\forwardvec}[1]{
	\settowidth\mystringlen{$#1$}
	\stackon[0.2ex]{#1}{
		\tikz{
			\draw[stroke=black,-{Latex[length=0.15em,width=0.2em]},line width=0.1ex] (0,0) -- (\mystringlen,0);
		}
	}
}
\newcommand{\backwardvec}[1]{
	\settowidth\mystringlen{$#1$}
	\stackon[0.2ex]{#1}{
		\tikz{
			\draw[stroke=black,{Latex[length=0.15em,width=0.2em]}-,line width=0.1ex] (0,0) -- (\mystringlen,0);
		}
	}
}
\newcommand{\forward}[1]{\forwardvec{#1}}
\newcommand{\backward}[1]{\backwardvec{#1}}
\newcommand{\RR}{\mathbb{R}}

\newcommand{\GG}{\mathcal{G}}

\newcommand{\NN}{\mathbb{N}}

\renewcommand{\thefootnote}{\fnsymbol{footnote}}
\renewcommand{\phi}{\varphi}
\renewcommand{\hat}{\widehat}

\theoremstyle{definition}
\newtheorem{thm}{Theorem}

\newtheorem{cor}[thm]{Corollary}

\newtheorem{claim}{Claim}
\crefname{claim}{claim}{claims}
\crefformat{claim}{#2Claim~#1#3}
\Crefformat{claim}{#2Claim~#1#3}
\crefname{obs}{Observation}{Observations}
\newtheorem*{lem*}{Lemma}
\theoremstyle{definition}
\newtheorem{conj}[thm]{Conjecture}
\newtheorem*{conj*}{Conjecture}

\def\lqedsymbol{\ifmmode$\lrcorner$\else{\unskip\nobreak\hfil
		\penalty50\hskip1em\null\nobreak\hfil$\rule{1ex}{1ex}$
		\parfillskip=0pt\finalhyphendemerits=0\endgraf}\fi} 
\newenvironment{claimproof}[1][\proofname]
{%
	\proof[#1]%
}
{%
	\endproof%
}

\presetkeys%
{todonotes}%
{inline}{}
\date{}

\begin{document}

\title{\bf\fontsize{18pt}{18pt}\selectfont 
Quasi-isometries between graphs with variable edge lengths}

\author{James Davies\,\footnotemark[1] \qquad Meike Hatzel\footnotemark[2] \qquad Robert Hickingbotham\,\footnotemark[3]}

\footnotetext[1]{University of Cambridge, Cambridge, UK (\texttt{jgd37@cam.ac.uk}).}

\footnotetext[2]{Discrete Mathematics Group, Institute for Basic Science (IBS), Daejeon, South Korea. E-mail: \href{research@meikehatzel.com}{research@meikehatzel.com}. Meike Hatzel's research was supported by the Institute for Basic Science (IBS-R029-C1).}

\footnotetext[3]{Univ. Lyon, ENS de Lyon, UCBL, CNRS, LIP, France (\texttt{robert.hickingbotham@ens-lyon.fr}).}

\maketitle

\begin{abstract}
This paper investigates quasi-isometries between graphs with variable edge lengths. A quasi-isometry is a mapping between metric spaces that approximately preserves distances, allowing for a bounded amount of additive and multiplicative distortion. Recently, Nguyen, Scott, and Seymour conjectured that, by appropriately adjusting the edge lengths of the target graph along with modifying the additive distortion constant, the multiplicative distortion factor could be eliminated. We disprove this conjecture.
\end{abstract}

% A quasi-isometry is a map between metric spaces that approximately preserves distances, allowing for a small amount of additive and multiplicative distortions. This paper investigates quasi-isometries between weighted graphs where distance is given by the shortest path length. Recently, Nguyen, Scott, and Seymour [2501.09839] conjectured that, up to a suitable change in the edge weighting of the target graph as well as change in the additive distortion, the multiplicative distortion could be eliminated. We disprove this conjecture by demonstrating that multiplicative distortion is, in fact, necessary for weighted graphs.

% Suggested shorter abstract:

% This paper investigates quasi-isometries between graphs with variable edge lengths. Recently, Nguyen, Scott, and Seymour conjectured that, up to a suitable change in the edge lengths of the target graph, the multiplicative distortion of the quasi-isometry could be negated at the cost of possibly increasing the additive distortion. We disprove this conjecture.

% Alternative shorter abstract:

\section{Introduction}
A quasi-isometry is a map between metric spaces that preserves the large-scale geometry of the spaces. Formally, given metric spaces $(X, d_X)$ and $(Y, d_Y)$, a map $\phi \colon X \to Y $ is an \defn{$(L, C)$-quasi-isometry} if there exist $L,C\in \NN$ such that, for all $ x_1, x_2 \in X $,  
\begin{equation*}
    L^{-1} \cdot d_X(x_1,x_2) - C \leq d_{Y}(\phi(x_1), \phi(x_2))\leq L \cdot d_X(x_1,x_2) + C,
\end{equation*}
and, for every $y\in Y$, there exists an $x \in X$ such that $d_{Y}(y,\phi(x)) \leq C$. If such a map exists, then we say that $(X, d_X)$ is \defn{$(L, C)$-quasi-isometric} to $(Y, d_Y)$.

Quasi-isometries play a central role in geometric group theory and metric geometry. They preserve large-scale geometric properties while ignoring small-scale differences. In particular, a large body of research in geometric group theory centres upon understanding which properties of groups are invariant under quasi-isometry.

This paper explores quasi-isometries between graphs (weighted and unweighted), where a \defn{weighted graph} is a pair $(H,w)$, where $H$ is a graph, and $w\colon E(H)\to \RR$ is an edge-weighting function. Distances are measured by the length of the shortest path where the length of a path is the sum of the weight of its edges. 

Recently, there has been significant interest in identifying conditions under which a graph is quasi-isometric to a simpler graph. This blossoming area, known as \emph{coarse graph theory}, seeks to use quasi-isometry to describe the large-scale geometry of complex graphs. 

A fundamental question concerning quasi-isometry is whether both the multiplicative distortion and additive distortion components are, in fact, necessary. First, consider whether the additive distortion factor $C$ is necessary (i.e. is it necessary to have $C\neq 0$). For general metric spaces, this component is indeed necessary. For example, the real line $\RR$ with the Euclidian distance is $(1,1)$-quasi-isometric to the integer lattice $\mathbb{Z}$, but not $(L,0)$-quasi-isometric to $\mathbb{Z}$ for any $L\in \NN$, since there will be a pair of arbitrarily close points in $\RR$ that would be mapped to distinct integers. 

However, when the metric spaces are induced by unweighted graphs, then the additive distortion factor can be avoided. Specifically, if a graph $G$ is $(L,C)$-quasi-isometric to another graph $H$, then it is $(L+1,0)$-quasi-isometric to $H$.

Now consider the distortion factor $L$ and whether we could always choose $L=1$. Such quasi-isometries are called \emph{slack-isometries} \cite{benjamini2013euclidean}. On one level, it is clear that we cannot always choose $L=1$. For instance, take a path $G$ of length $C+1$ and another path $H$ of length $2C+2$.
Then $G$ is $(2,0)$-quasi-isometric to $H$ but is not $(1,C)$-quasi-isometric to $H$.
However, the situation changes drastically when we allow $H$ to come from a family of graphs.
For example, a remarkable theorem by \citet{CDNRV2012} and also \citet{Kerr2023} states that for the class of trees, additive distortion alone suffices.

\begin{thm}[\cite{CDNRV2012,Kerr2023}]
	\label{thm:tree_case}
    For all $L,C\in \NN$, there exists $C'\in \NN$ such that if a graph $G$ is $(L, C)$-quasi-isometric to a tree, then $G$ is $(1,C')$-quasi-isometric to a tree.
\end{thm}

More recently, \citet{NSS2025treewidth} proved a similar result for graphs of bounded pathwidth in terms of weighted graphs. 

\begin{thm}[\cite{NSS2025treewidth}]
	\label{thm:pathwidth}
    For all $L,C,k\in \NN$, there exists $C'\in \NN$ such that if $\phi$ is an $(L, C)$-quasi-isometry from a graph $G$ to a graph $H$ with pathwidth at most $k$, then there is an edge-weighting function $w\colon E(H) \to \mathbb{Z}$ such that the same function $\phi$ is a $(1, C')$-quasi-isometry from $G$ to the weighted graph $(H,w)$.
\end{thm}

In light of this result, \citet{NSS2025treewidth} conjectured that, up to a suitable edge weighting function, additive distortion alone suffices. Specifically, they conjectured the following: 

\begin{conj}[\cite{NSS2025treewidth}]\label{MainConjecture}
    For all $L,C\in \NN$, there exists $C'\in \NN$ such that if $\phi$ is an $(L, C)$-quasi-isometry from a graph $G$ to a graph $H$, then there is an edge-weighting function $w\colon E(H) \to \mathbb{N}$ such that the same function $\phi$ is a $(1, C')$-quasi-isometry from $G$ to the weighted graph $(H,w)$.
\end{conj}

This conjecture, if true, would imply that for any class of graphs $\mathcal{H}$ closed under contracting edges and taking subdivisions, if a graph $G$ is $(L,C)$-quasi-isometric to a graph $H\in \mathcal{H}$, then $G$ is $(1,C')$-quasi-isometric to a graph $H'\in \mathcal{H}$ where $C'$ depends only on $L$ and $C$. This would be a remarkably powerful result if true.

In this paper, we disprove this conjecture by constructing explicit counterexamples demonstrating that multiplicative distortion for weighted graphs is, in fact, necessary.

\begin{thm}\label{MainTheorem}
     For every $C\in \mathbb{N}$, there exist graphs $G$ and $H$ and a $(2,1)$-quasi-isometry $\phi\colon V(G)\to V(H)$ such that, for every edge weighting $w\colon E(H)\to \mathbb{R}^+$ of $H$, the map $\phi$ is not a $(1,C)$-quasi-isometry from $G$ to $(H,w)$.
\end{thm}

Our proof for \cref{MainTheorem} is based on orientated graphs $H$ with large girth and chromatic number. We use the orientation of the edges to split each vertex into a new edge to define our graph $G$ (see \cref{fig:construction}). Since $H$ can be obtained from $G$ by contracting disjoint edges, the natural map $\phi$ between $G$ and $H$ is a $(2,1)$-quasi-isometry. Next, we assume that $H$ is given an edge-weighting function $w\colon E(H)\to \RR^+$. We then separate the light-weight edges from the heavy-weight edges and then find different types of long orientated paths within the graph. This allows us to find in $H$ either a long path of light-weight edges that traverse many new edges, or a long path of heavy-weight edges that avoids the new edges. Since $H$ has large girth, such paths are geodesic which allows us to contradict $\phi$ being a $(1,C)$-quasi-isometry from $G$ to $(H,w)$.

\subsection{Preliminaries}

Let $G$ be a graph. The \defn{girth} of $G$ is the length of a shortest cycle in $G$. For $k\in \NN$, a \defn{proper $k$-colouring} of a $G$ is a function $c \colon V(G)\to \{1,\dots,k\}$ such that $c(u)\neq c(v)$ whenever $uv\in E(G)$. The \defn{chromatic number} $\chi(G)$ is the minimum $k\in \NN$ for which $G$ has a proper $k$-colouring.

Let $(H,w)$ be a weighted graph. For a path $P=(v_0,v_1,\dots,v_n)$ in $(H,w)$, we say that the  \defn{length} of $P$ is $\sum_{i=1}^n w(v_{i-1}v_i)$ and the \defn{hop-length} of $P$ is $n$. The path $P$ is \defn{geodesic} if it is a path of minimum length in $(H,w)$ between $v_0$ and $v_n$.

An \defn{oriented graph} $\forward{H\hspace*{1pt}}$ is a graph where each edge has a direction. The \defn{chromatic number} of $\forward{H\hspace*{1pt}}$ is the chromatic number of $H$. Let $\forward{P_n}=(v_0,v_1,v_2,\dots,v_n)$ be an oriented path. We say that $\forward{P_n}$ is \defn{directed} if $\forward{v_{i-1}v_{i}}$ for every $i\in [n]$ and we say that $\forward{P_n}$ is \defn{alternating} if $\forward{v_{i-1}v_{i}}$ for every odd $i\in [n]$ and $\backward{v_{i-1}v_{i}}$ for every even $i\in [n]$.

\section{Proof}

The following classical result of \citet{burr1980subtrees} allows us to find long oriented paths in graphs with large chromatic number.

\begin{thm}[\cite{burr1980subtrees}]\label{OrientedPathLemma}
     For every $n\in \mathbb{N}$, for every oriented path $\forward{P_n}$ on $n$ vertices, every oriented graph $\forward{G\hspace*{1pt}}$ that does not contain $\forward{P_n}$ as an oriented subgraph has chromatic number at most $n^2$.
\end{thm}

We now prove our main result.

\begin{proof}[Proof of \cref{MainTheorem}]
Let $H$ be a graph with girth at least $40C^2$ and chromatic number at least $256C^4$. By a classical result of \citet{Erdos59}, such a graph exists. Define an oriented graph $\forward{H\hspace*{1pt}}$ by arbitrarily orientating each edge of $H$. Construct a new graph $G$ from $H$ by replacing each vertex $v \in V(H)$ with two adjacent vertices $v^{-}$ and $v^{+}$, and adding the edge $u^{+} v^{-}$ whenever $\forward{uv} \in E(\forward{H\hspace*{1pt}})$.
See \cref{fig:construction} for an illustration of this construction.
Since $H$ has girth at least $40C^2$, so does $G$.

\begin{figure}[!b]
	\definecolor{vcolor}{rgb}{0.68, 0.05, 0.0}
	\begin{tikzpicture}
		\path[boundingBox] (-7.8,-2.8) rectangle (7.8,2.3);
		\node (H-pic) at (-5.5,0) {
			\begin{tikzpicture}[scale=1.2]
				\node[vertex,vcolor,inner sep=.8mm] (v) at (0,0) {};
				\node (v-label) at ($(v)+(45:0.3)$) {\textcolor{vcolor}{$v$}};
				\node[vertex] (x-1) at ($(v)+(170:1)$) {};
				\node[vertex] (x-2) at ($(v)+(95:1)$) {};
				\node[vertex] (x-3) at ($(v)+(350:1)$) {};
				\node[vertex] (x-4) at ($(v)+(275:1)$) {};
				\node[vertex] (x-5) at ($(v)+(215:1)$) {};
				\node[vertex] (x-p-5) at ($(x-5)+(175:0.6)$) {};
				\node (y-1) at ($(v)+(165:2)$) {};
				\node (y-2) at ($(v)+(185:2)$) {};
%				\node (y-3) at ($(v)+(185:2)$) {};
				\node (y-4) at ($(v)+(260:2)$) {};
				\node (y-5) at ($(v)+(280:2)$) {};
				\node (y-6) at ($(v)+(330:2)$) {};
				\node (y-7) at ($(v)+(350:2)$) {};
				\node (y-8) at ($(v)+(370:2)$) {};
				\node (y-9) at ($(v)+(70:2)$) {};
				\node (y-10) at ($(v)+(100:2)$) {};
				
				\draw[edge] (v) to (x-1);
				\draw[edge] (v) to (x-2);
				\draw[edge] (v) to (x-3);
				\draw[edge] (v) to (x-4);
				\draw[edge] (v) to (x-5);
				\draw[edge] (x-5) to (x-p-5);
				\draw[edge] (x-1) to (y-1.center);
				\draw[edge] (x-1) to (y-2.center);
				\draw[edge] (x-4) to (y-4.center);
				\draw[edge] (x-4) to (y-5.center);
				\draw[edge] (x-3) to (y-6.center);
				\draw[edge] (x-3) to (y-7.center);
				\draw[edge] (x-3) to (y-8.center);
				\draw[edge] (x-2) to (y-9.center);
				\draw[edge] (x-2) to (y-10.center);
				
				\draw[fill=none,dashed] (v) circle (2.0);
			\end{tikzpicture}};
		\node (H-pic-label) at (-4.5,-2.7) {$H\vphantom{\forward{H}}$};
		\node (dirH-pic) at (0,0) {
			\begin{tikzpicture}[scale=1.2]
				\node[vertex,vcolor,inner sep=.8mm] (v) at (0,0) {};
				\node (v-label) at ($(v)+(45:0.3)$) {\textcolor{vcolor}{$v$}};
				\node[vertex] (x-1) at ($(v)+(170:1)$) {};
				\node[vertex] (x-2) at ($(v)+(95:1)$) {};
				\node[vertex] (x-3) at ($(v)+(350:1)$) {};
				\node[vertex] (x-4) at ($(v)+(275:1)$) {};
				\node[vertex] (x-5) at ($(v)+(215:1)$) {};
				\node[vertex] (x-p-5) at ($(x-5)+(175:0.6)$) {};
				\node (y-1) at ($(v)+(165:2)$) {};
				\node (y-2) at ($(v)+(185:2)$) {};
				%				\node (y-3) at ($(v)+(185:2)$) {};
				\node (y-4) at ($(v)+(260:2)$) {};
				\node (y-5) at ($(v)+(280:2)$) {};
				\node (y-6) at ($(v)+(330:2)$) {};
				\node (y-7) at ($(v)+(350:2)$) {};
				\node (y-8) at ($(v)+(370:2)$) {};
				\node (y-9) at ($(v)+(70:2)$) {};
				\node (y-10) at ($(v)+(100:2)$) {};
				
				\draw[directededge] (v) to (x-2);
				\draw[directededge] (v) to (x-3);
				\draw[directededge] (x-1) to (v);
				\draw[directededge] (x-4) to (v);
				\draw[directededge] (v) to (x-5);
				\draw[directededge] (x-p-5) to (x-5);
				\draw[directededge] (x-1) to (y-1.center);
				\draw[directededge] (y-2.center) to (x-1);
				\draw[directededge] (x-4) to (y-4.center);
				\draw[directededge] (y-5.center) to (x-4);
				\draw[directededge] (x-3) to (y-6.center);
				\draw[directededge] (y-7.center) to (x-3);
				\draw[directededge] (x-3) to (y-8.center);
				\draw[directededge] (y-9.center) to (x-2);
				\draw[directededge] (y-10.center) to (x-2);
				
				\draw[fill=none,dashed] (v) circle (2.0);
		\end{tikzpicture}};
		\node (dirH-pic-label) at ($(H-pic-label)+(5.5,0)$) {$\forward{H}$};
		\node (G-pic) at (5.5,0) {
			\begin{tikzpicture}[scale=1.2]
				\node (v) at (0,0) {};
				\node[vertex,vcolor,inner sep=.6mm] (v-in) at ($(v)+(215:0.3)$) {};
				\node (v-in-label) at ($(v-in)+(200:0.3)$) {\textcolor{vcolor}{$v^{-}$}};
				\node[vcolor,vertexOut,inner sep=.6mm] (v-out) at ($(v)+(45:0.3)$) {};
				\node (v-out-label) at ($(v-out)+(40:0.35)$) {\textcolor{vcolor}{$v^{+}$}};
%				\node (v-label) at ($(v)+(45:0.3)$) {\textcolor{vcolor}{$v$}};

				\node (x-1) at ($(v)+(170:1)$) {};
				\node[vertex] (x-1-in) at ($(x-1)+(180:0.25)$) {};
				\node[vertexOut] (x-1-out) at ($(x-1)+(0:0.25)$) {};
				
				\node (x-2) at ($(v)+(95:1)$) {};
				\node[vertex] (x-2-in) at ($(x-2)+(10:0.25)$) {};
				\node[vertexOut] (x-2-out) at ($(x-2)+(190:0.25)$) {};
				
				\node (x-3) at ($(v)+(350:1)$) {};
				\node[vertex] (x-3-in) at ($(x-3)+(100:0.25)$) {};
				\node[vertexOut] (x-3-out) at ($(x-3)+(280:0.25)$) {};
				
				\node (x-4) at ($(v)+(275:1)$) {};
				\node[vertex] (x-4-in) at ($(x-4)+(-20:0.25)$) {};
				\node[vertexOut] (x-4-out) at ($(x-4)+(160:0.25)$) {};
				
				\node (x-5) at ($(v)+(215:1)$) {};
				\node[vertex] (x-5-in) at ($(x-5)+(80:0)$) {};
				\node[vertexOut] (x-5-out) at ($(x-5)+(260:0.4)$) {};
				
				\node (x-p-5) at ($(x-5)+(175:0.6)$) {};
				\node[vertex] (x-p-5-in) at ($(x-p-5)+(220:0.2)$) {};
				\node[vertexOut] (x-p-5-out) at ($(x-p-5)+(60:0.2)$) {};
				
				\node (y-1) at ($(v)+(165:2)$) {};
				\node (y-2) at ($(v)+(185:2)$) {};
				%				\node (y-3) at ($(v)+(185:2)$) {};
				\node (y-4) at ($(v)+(260:2)$) {};
				\node (y-5) at ($(v)+(280:2)$) {};
				\node (y-6) at ($(v)+(330:2)$) {};
				\node (y-7) at ($(v)+(350:2)$) {};
				\node (y-8) at ($(v)+(370:2)$) {};
				\node (y-9) at ($(v)+(70:2)$) {};
				\node (y-10) at ($(v)+(100:2)$) {};
				
				\draw[edge,vcolor] (v-in) to (v-out);
				\draw[edge] (x-1-in) to (x-1-out);
				\draw[edge] (x-2-in) to (x-2-out);
				\draw[edge] (x-3-in) to (x-3-out);
				\draw[edge] (x-4-in) to (x-4-out);
				\draw[edge] (x-5-in) to (x-5-out);
				\draw[edge] (x-p-5-in) to (x-p-5-out);
				
				\draw[edge] (v-out) to (x-2-in);
				\draw[edge] (v-out) to (x-3-in);
				\draw[edge] (x-1-out) to (v-in);
				\draw[edge] (x-4-out) to (v-in);
				\draw[edge,out=300,in=-30] (v-out) to (x-5-in);
				\draw[edge] (x-p-5-out) to (x-5-in);
				\draw[edge,out=90,in=350] (x-1-out) to (y-1.center);
				\draw[edge] (y-2.center) to (x-1-in);
				\draw[edge] (x-4-out) to (y-4.center);
				\draw[edge] (y-5.center) to (x-4-in);
				\draw[edge] (x-3-out) to (y-6.center);
				\draw[edge] (y-7.center) to (x-3-in);
				\draw[edge,out=10,in=200] (x-3-out) to (y-8.center);
				\draw[edge] (y-9.center) to (x-2-in);
				\draw[edge] (y-10.center) to (x-2-in);
				
				\draw[fill=none,dashed] (v) circle (2.0);
		\end{tikzpicture}};
		\node (G-pic-label) at ($(dirH-pic-label)+(5.5,0)$) {$G\vphantom{\forward{H}}$};
	\end{tikzpicture}
	\caption{Starting with $H$, we orientate the edges to obtain \(\protect\forward{H}\) then we split the vertices to obtain $G$.}
	\label{fig:construction}
\end{figure}
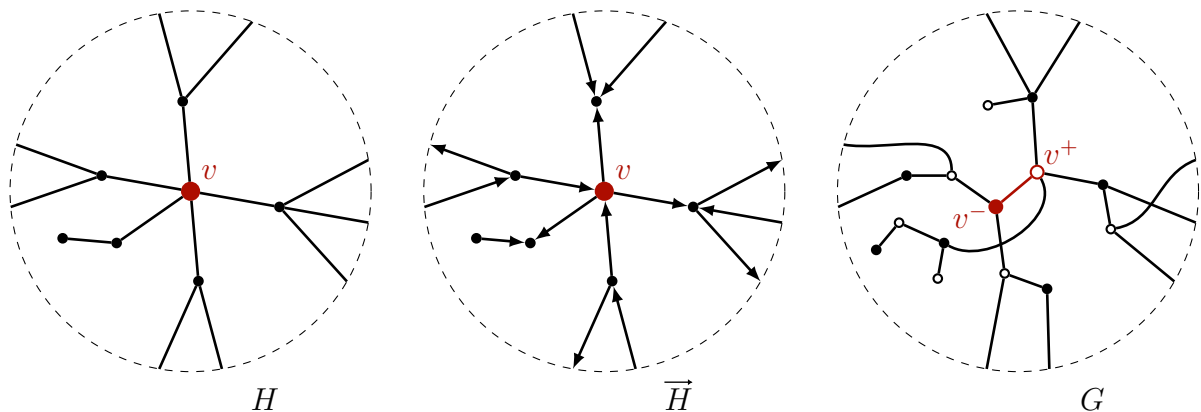

\begin{claim}
	\label{claim:short-length-implies-geodesic}
	Every path in $H$ (or in $G$) of length at most $10C$ is geodesic.
\end{claim}
\begin{claimproof}
    Let $P$ be a path in $H$ of length at most $10C$. Suppose $H$ contains a shorter path $P'$ with the same end-vertices as $P$. Then $P \cup P'$ would contain a cycle of length at most $20C$, contradicting $H$ having girth at least $40C^2$. The same argument also applies if $P$ was a path in $G$.
\end{claimproof}

Define $\phi\colon V(G)\to V(H)$ by mapping $v^{+}$ and $v^{-}$ to $v$ for each $v\in V(H)$. Since $H$ is obtained from $G$ by contracting vertex-disjoint edges, it follows that for all $x,y\in V(G)$, we have:
$$\dist_H(\phi(x),\phi(y))\leq \dist_G(x,y)\leq 2 \cdot \dist_H(\phi(x),\phi(y))+1.$$ 
Thus, $\phi$ is a $(2,1)$-quasi-isometry from $G$ to $H$.

Assume, for the sake of contradiction, that there exists an edge-weighting function $w\colon E(H)\to \mathbb{R}$ such that $\phi\colon V(G)\to V(H)$ is a $(1,C)$-quasi-isometry from $G$ to $(H,w)$.

\begin{claim}
	\label{claim:weight}
	Every edge in $(H,w)$ has weight at most $C+1$.
\end{claim}
\begin{claimproof}
     Suppose there exists an edge $uv\in E(H)$ with $w(uv)> C+1$.
     Since either $u^+v^-$ or $u^-v^+$ is in $G$ and $\phi$ is a $(1,C)$-quasi-isometry to $(H,w)$, we have $\dist_{(H,w)}(u,v)\leq C+1$.
     So $(H,w)$ contains a $(u,v)$-path $P_{uv}$ with length at most $C+1$.
     This means that for every $x\in V(P_{uv})$, we have $\dist_{(H,w)}(u,x)\leq C+1$.
     Given $H$ has girth at least $40C^2$, the path $P_{uv}$ has hop-length at least $40C^2-1$. Consequently, there is a vertex $x\in V(P_{uv})$ such that $\dist_H(u,x)\geq 2C+2$. Therefore, $\dist_{G}(u^+,x^+)\geq \dist_H(u,x)\geq 2C+2$, while simultaneously $\dist_{(H,w)}(\phi(u^+),\phi(x^+))\leq C+1$, which contradicts the assumption that $\phi$ is a $(1,C)$-quasi-isometry. 
\end{claimproof}

\begin{claim}
	\label{claim:short-hop-length-implies-geodesic}
	Every path in $(H,w)$ of hop-length at most $4C$ is geodesic.
\end{claim}
\begin{claimproof}
    Suppose, for the sake of contradiction, that there exists a path $(u,v)$-path $P$ in $(H,w)$ of hop-length at most $4C$ that is not geodesic. Assume $P$ is of minimal hop-length. By \cref{claim:weight}, we have $\dist_{(H,w)}(u,v)\leq 4C(C+1)\leq 8C^2$. Since $P$ is not geodesic, there exists an alternative $(u,v)$-path $P'$ with length less than $8C^2$. This means that for every $x\in V(P_{uv})$, we have $\dist_{(H,w)}(u,x)\leq 8C^2$. By the minimality of $P$, $P'$ is internally disjoint from $P$, so $P\cup P'$ is a cycle. Given that $H$ has girth at least $40C^2$, $P'$ has hop-length at least $30C^2$. Thus, there exists a vertex $x\in V(P')$ such that $\dist_H(u,x)\geq 10C^2$, but $\dist_{(H,w)}(u,x)\leq 8C^2$. This leads to a contradiction in the same manner as before, proving the claim.
\end{claimproof}

Define an edge $e\in E(H)$ as \emph{light} if $w(e)\leq 1.5$ and \emph{heavy} otherwise. Let $X_1$ and $X_2$ be the spanning subgraphs of $H$ consisting of light and heavy edges, respectively. Since $E(H) = E(X_1) \cup E(X_2)$, we have:
$$\max\{\chi(X_1),\chi(X_2)\}\geq 16C^2.$$
Otherwise, we could $(16C^2-1)$-colour each subgraph separately and obtain a proper $(16C^2-1)^2$-colouring of $H$, contradicting $\chi(H) \geq 256C^4$.

Choose $i\in \{1,2\}$ such that $\chi(X_i)\geq 16C^2$. Let $\forward{X}_i$ be the oriented subgraph of $\forward{H\hspace*{1pt}}$ corresponding to $X_i$. By \cref{OrientedPathLemma}, $\forward{X}_i$ contains both a directed path $P=(v_0,v_1,\dots,v_{4C})$ of hop-length $4C$ and an alternating path $\hat{P}=(u_0,u_1,\dots,u_{4C})$ of hop-length $4C$.

If $i=1$ (all edges are light), then $\dist_{(H,w)}(v_0,v_{4C})\leq 1.5\cdot 4C=6C$. Since $P$ is a directed path, $G$ contains a $(v_0^{-},v_{4C}^+)$-path with length $8C$ of the form:
$$(v_0^-,v_0^+,v_1^-,v_1^+,\dots,v_{4C}^-,v_{4C}^+).$$
By \cref{claim:short-length-implies-geodesic}, this path is geodesic in $G$.
This contradicts $\phi$ being a $(1,C)$-quasi-isometry since $$\dist_G(v_0^{-},v_{4C}^+)=8C > 6C \geq \dist_{(H,w)}(\phi(v_0^{-}),\phi(v_{4C}^+))+C.$$

If $i=2$ (all edges are heavy), then by \cref{claim:short-hop-length-implies-geodesic}, the path $\hat{P}$ is geodesic in $(H,w)$. Therefore, $\dist_{(H,w)}(\hat{x},\hat{y})\geq 1.5\cdot 4C= 6C$. Since $\hat{P}$ is an alternating path in $\forward{H\hspace*{1pt}}$, $G$ contains a $(u_0^{+},u_{4C}^+)$-path with length $4C$ of the form:
$$(u_0^+,u_1^-,u_2^+,\dots,u_{4C-1}^-,u_{4C}^+).$$
This again contradicts $\phi$ being a $(1,C)$-quasi-isometry since 
$$\dist_G(u_0^{+},u_{4C}^+)=4C < 6C \leq \dist_{(H,w)}(\phi(u_0^{+}),\phi(u_{4C}^+))-C.$$ 
This completes the proof.
\end{proof}

\section{Conclusion}
We conclude by presenting two modifications to the graph $G$ and $H$ in \cref{MainTheorem} that further strengthen the counterexample. 

First, the condition of $G$ being `$(2,1)$-quasi-isometric' can be refined to `$(1+\epsilon,1)$-quasi-isometric' for any given $\epsilon>0$. This can be achieved by replacing each edge of $H$ with a path of length $\lceil 1/\epsilon \rceil$ before splitting the original vertices of $H$ to define $G$. The rest of the argument in \cref{MainTheorem} then goes through with some minor adjustments.

\begin{cor}
	\label{cor:1plus-epsilon-version}
	For every $C\in \mathbb{N}$ and $\epsilon>0$, there exist graphs $G$ and $H$ and a $(1+\epsilon,1)$-quasi-isometry $\phi$ from $G$ to $H$ such that, for every edge weighting $w\colon E(H)\to \mathbb{R}$ of $H$, the map $\phi$ is not a $(1,C)$-quasi-isometry from $G$ to $(H,w)$.
\end{cor}

Second, we can ensure that no map between $G$ and $(H,w)$ qualifies as a $(1,C)$-quasi-isometry by attaching long paths of significantly different lengths to each vertex of $H$ before splitting the original vertices to define $G$. This modification forces any potential $(1,C)$-quasi-isometry $\phi\colon V(G)\to V(H)$ to map $\phi(v^+)$ and $\phi(v^-)$ close to $v$ for each $v\in V(H)$. By continuing with the argument in \cref{MainTheorem}, we conclude that no edge weighting of $H$ allows $\phi$ to be a $(1,C)$-quasi-isometry.

\begin{cor}
	\label{cor:any-map-version}
	For every $C\in \mathbb{N}$ and $\epsilon>0$, there exist graphs $G$ and $H$ such that $G$ is $(1+\epsilon,1)$-quasi-isometric to $H$ but there is no edge weighting $w\colon E(H)\to \mathbb{R}$ of $H$ such that $G$ is $(1,C)$-quasi-isometric to $(H,w)$.
\end{cor}

Finally, while \cref{MainConjecture} is false, we would like to point out the following weakening of the conjecture which remains open. 

\begin{conj}[\cite{NSS2025treewidth}]\label{ConjectureGraphClasses}
    Let $\GG$ be a class of connected graphs closed under contracting edges and subdividing edges. For all $L,C\in \NN$, there exists $C'\in \NN$ such that if a graph $G$ is $(L, C)$-quasi-isometric to a graph in $\GG$, then $G$ is $(1, C')$-quasi-isometric to a graph in $\GG$.
\end{conj}

A special case of \cref{ConjectureGraphClasses} that is of particular interest is when $\GG$ is the class of connected planar graphs, as possibility first conjectured by Georgakopoulos (private communication).

% \begin{thm}
%     For every $C\in \mathbb{N}$ and $\epsilon>0$, there exists a graph $G$ that is $(1+\epsilon,1)$-quasi-isometric to a graph $H''$ such that, for every edge weighting $w\colon E(H)\to \mathbb{R}$ of $H$, every map $\phi\colon V(G)\to V(H)$ is not a $(1,C)$-quasi-isometry from $G$ to $(H,w)$.
% \end{thm}
% \begin{proof}[Proof Sketch]
%     The proof follows the same strategy as that in \cref{MainTheorem}. We highlight the key difference.
%     First, we explain how to construct $G$ so that it is $(1+\epsilon,1)$-quasi-isometric to a graph $H'$. As in the proof of \cref{MainTheorem}, we choose $H$ to be a graph with sufficiently large girth and chromatic number.
    
%     Before splitting the vertices of $H$, we first replace each edge of $H$ by a path of length $\ceil{10/\epsilon}$ to obtain the graph $H'$. The graph $G$ is then obtained from $H'$ by orientating the original edges of $H$ and then splitting each vertex in $V(H)$. By doing so, the graph $G$ is then $(1+\epsilon,1)$-quasi-isometric to $H'$.

%     Now to ensure that every map is not a $(1,C)$-quasi-isometry, we then modifying the graph $H'$ further by attaching long paths of significantly distinct length to each vertex of $H'$. In doing so, this forces any $(L,C)$-quasi-isometry from $G$ to $(H,w)$ to be essentially the same; that is, if $\phi_1$ and $\phi_2$ are $(1,C)$-quasi-isometry from $G$ to $(H,w)$, then $\dist_{(H,w)}(\phi_1(v),\phi_2(v))\leq f(C)$ for all $v\in V(G)$.  The rest of the proof then follows as in \cref{MainTheorem}.
% \end{proof}

\paragraph{Acknowledgement:} This work was completed at the 12th Annual Workshop on Geometry and Graphs held at Bellairs Research Institute in February 2025. Thanks to the organisers and participants for providing a stimulating work environment.

\renewcommand{\thefootnote}{\arabic{footnote}}

{
\fontsize{11pt}{12pt}
\selectfont
	
\hypersetup{linkcolor={red!70!black}}
\setlength{\parskip}{2pt plus 0.3ex minus 0.3ex}

\bibliographystyle{DavidNatbibStyle}
\bibliography{RobReferences}

\begin{thebibliography}{6}
\providecommand{\natexlab}[1]{#1}
\providecommand{\msn}[1]{MR:\,\href{http://www.ams.org/mathscinet-getitem?mr=MR{#1}}{#1}}
\providecommand{\ZBL}[1]{Zbl:\,\href{https://www.zentralblatt-math.org/zmath/en/search/?q=an:#1}{#1}}
\providecommand{\url}[1]{\texttt{#1}}
\providecommand{\urlprefix}{}
\expandafter\ifx\csname urlstyle\endcsname\relax
  \providecommand{\doi}[1]{doi:\discretionary{}{}{}#1}\else
  \providecommand{\doi}{doi:\discretionary{}{}{}\begingroup
  \urlstyle{rm}\Url}\fi

\bibitem[{Benjamini(2013)}]{benjamini2013euclidean}
\textsc{Itai Benjamini}.
\newblock
  \href{https://link.springer.com/chapter/10.1007/978-3-642-39286-3_2}{Euclidean
  vs. graph metric}.
\newblock In \emph{Erdős Centennial}, vol.~25 of \emph{Bolyai Soc. Math.
  Stud.}, pp. 35--57. Springer, 2013.

\bibitem[{Burr(1980)}]{burr1980subtrees}
\textsc{Stefan~A. Burr}.
\newblock Subtrees of directed graphs and hypergraphs.
\newblock In \emph{Proceedings of the Eleventh Southeastern Conference on
  Combinatorics, Graph Theory and Computing}, vol.~28 of \emph{Congressus
  Numerantium}, pp. 227--239. 1980.

\bibitem[{Chepoi et~al.(2012)Chepoi, Dragan, Newman, Rabinovich, and
  Vax{\`e}s}]{CDNRV2012}
\textsc{Victor Chepoi, Feodor~F. Dragan, Ilan Newman, Yuri Rabinovich, and Yann
  Vax{\`e}s}.
\newblock \href{https://doi.org/10.1007/s00454-011-9386-0}{Constant
  approximation algorithms for embedding graph metrics into trees and
  outerplanar graphs}.
\newblock \emph{Discrete \& Computational Geometry}, 47:187--214, 2012.

\bibitem[{Erd{\H{o}}s(1959)}]{Erdos59}
\textsc{Paul Erd{\H{o}}s}.
\newblock \href{https://doi.org/10.4153/CJM-1959-003-9}{Graph theory and
  probability}.
\newblock \emph{Canad. J. Math.}, 11:34--38, 1959.

\bibitem[{Kerr(2023)}]{Kerr2023}
\textsc{A.~Kerr}.
\newblock \href{https://doi.org/10.4171/GGD/733}{Tree approximation in
  quasi-trees}.
\newblock \emph{Groups, Geometry and Dynamics}, 17:1193--1233, 2023.

\bibitem[{Nguyen et~al.(2025)Nguyen, Scott, and Seymour}]{NSS2025treewidth}
\textsc{Tung Nguyen, Alex Scott, and Paul Seymour}.
\newblock \href{http://arxiv.org/abs/2501.09839}{Coarse tree-width}.
\newblock 2025, arXiv:2501.09839.

\end{thebibliography}
}

\end{document}